\documentclass[11pt]{article}

\usepackage{amsthm,amsmath,hyperref,txfonts,mathrsfs}
\usepackage{multirow}
\usepackage{color}
\usepackage{authblk}
\usepackage[marginal]{footmisc}

\newtheorem{theorem}{Theorem}
\newtheorem{lemma}[theorem]{Lemma}

\newtheorem{remark}[theorem]{Remark}

\allowdisplaybreaks

\begin{document}
\title{Corrigendum to ``Optimal time-decay estimates for an Oldroyd-B model
with zero viscosity [J. Differential Equations, 306(2022), 456--491]''}

\author[a]{Jinrui Huang}
\author[b]{Yinghui Wang}
\author[b]{Huanyao Wen\thanks{Corresponding author.\\ $\quad\quad\quad$ E-mail addresses: huangjinrui1@163.com (J. Huang), yhwangmath@163.com (Y. Wang), mahywen@scut.edu.cn (H. Wen), rzz@ccnu.edu.cn (R. Zi).}}
\author[c]{Ruizhao Zi}
\affil[a]{School of Mathematics and Computational Science, Wuyi University, Jiangmen 529020, China. }
\affil[b]{School of Mathematics, South China University of Technology, Guangzhou 510641, China.}
\affil[c]{School of Mathematics and Statistics $\&$ Hubei Key Laboratory of Mathematical  Sciences,  Central China Normal University, Wuhan 430079, China.}

\date{}
\maketitle
\begin{abstract}
	The present note is to make minor correction  on the assumption of Theorem 1.2 and its proof in our paper [Jinrui Huang, Yinghui Wang,  Huanyao Wen and Rizhao Zi, {\it J. Differential Equations}, 306(2022), 456--491].
\end{abstract}




\setcounter{section}{0} \setcounter{equation}{0}

\section{Reasons to make corrections}
In the original publication of paper \cite{HWWZ_2022_JDE}, the authors studied the optimal time-decay estimates for the following Cauchy problem of Oldroyd-B system
\begin{eqnarray} \label{system}
\begin{cases}
     \partial_tu+u\cdot\nabla u+\nabla p-\epsilon\Delta u=\kappa {\rm div}\tau,\\
     \partial_t\tau+u\cdot\nabla\tau-\mu\Delta\tau+\beta\tau=Q(\nabla u,\tau)+\alpha\mathbb{D}u,\\
     {\rm div}u=0,\\
     (u,\tau)(x,0)=(u_0,\tau_0),
\end{cases}
\end{eqnarray} on $\mathbb{R}^3\times (0,\infty)$. Here, using the same notations as in \cite{HWWZ_2022_JDE}, $u=(u_1,u_2,u_3)^\top:\mathbb{R}^3\rightarrow\mathbb{R}^3$ the velocity field of fluid, $\tau\in \mathbb{S}_3(\mathbb{R})$ is the purely elastic (the polymer) part of the stress tensor, $p\in\mathbb{R}$ is the pressure function of the fluid, $\mathbb{D}u=\frac12\left(\nabla u+\left(\nabla u\right)^\top\right)$ is the deformation tensor, and $$Q(\nabla u,\tau)=\Omega\tau-\tau\Omega+b(\mathbb{D}u \tau+\tau\mathbb{D}u),$$ admits the usual bilinear form  with the skew-symmetric part of velocity gradient $\Omega=\frac12\left(\nabla u-\left(\nabla u\right)^\top\right)$ and some $b\in [-1,1]$. The parameters $\epsilon$, $\kappa$, $\mu$, $\beta$ and $\alpha$ satisfy that $\epsilon, \mu\geq0, \kappa,\beta,\alpha>0$.  The following two cases were studied
\begin{itemize}
	\item Case I: $\mu>0$ and $\epsilon\ge0$;
	
	\item Case II: $\epsilon>0$ and $\mu\ge0$.
\end{itemize}
In pages 460--461 of paper \cite{HWWZ_2022_JDE}, the authors stated the following results,

\begin{theorem}[Theorem 1.2 of  \cite{HWWZ_2022_JDE}]
	\label{theorem_decay}
	For any given $\epsilon$ and $\mu$ satisfying the Case I or Case II, letting $(u^{\epsilon,\mu},\tau^{\epsilon,\mu})$ be the solution as in Theorem 1.1, then the following optimal time-decay estimates hold.

	\medskip

	\noindent {\em (i)} Assume that $(u_0,\tau_0)\in L^1(\mathbb{R}^3)$. Then we have upper time-decay estimates of the solution as below:
	\begin{equation}\label{opti1}
	\ \|\nabla^ku^{\epsilon,\mu}(t)\|_{L^2}\leq C_2(1+t)^{-\frac34-\frac{k}{2}},\ k=0,1,2,
	\tag{1.7}
	\end{equation}
	and
	\begin{equation}\label{opti2}
	\ \|\nabla^{k_1}\tau^{\epsilon,\mu}(t)\|_{L^2}\leq C_2(1+t)^{-\frac54-\frac{k_1}{2}},\ k_1=0,1\tag{1.8}
	\end{equation}
	for any $t\geq0$ and generic positive constant $C_2$ which depends only on $\|(u_0,\tau_0)\|_{H^3\cap L^1}$ and $C_1$.
	
\medskip

	\noindent {\em(ii)} Assume that $(u_0,\tau_0)\in L^1(\mathbb{R}^3)$ and in addition that $\inf_{0\leq |\xi|\leq R}|\hat{u}_0|\geq c_{0}>0$. Then there exist a positive time $t_0=t_0(\alpha,\kappa,\beta,\|(u_0,\tau_0)\|_{L^1})$ and a positive generic constant $c=c(\alpha,\kappa,\beta,c_0,C_2)$, such that
	\begin{equation}\label{opti3}
	\|\nabla^ku^{\epsilon,\mu}(t)\|_{L^2}\geq c(1+t)^{-\frac34-\frac{k}{2}},\ k=0,1,2,\tag{1.9}
	\end{equation}
	and
	\begin{equation}\label{opti4}
	\|\nabla^{k_1}\tau^{\epsilon,\mu}(t)\|_{L^2}\geq c(1+t)^{-\frac54-\frac{k_1}{2}},\ k_1=0,1 \tag{1.10}
	\end{equation}
	for any $t\geq t_0$.
	
\end{theorem}
The authors apologize to inform that the assumption $\inf_{0\leq |\xi|\leq R}|\hat{u}_0|\geq c_{0}>0$ in Theorem 1.2 (ii) is contradictory with the assumptions $u_0\in H^3(\mathbb{R}^3)\cap L^1(\mathbb{R}^3)$ and $\mathrm{div} u_0 = 0$. The reason will be stated as follows.

As pointed out in \cite{Schonbek_Suli_1996}, the following Lemma holds.

\begin{lemma}[Lemma B in \cite{Schonbek_Suli_1996}, (page 721)] \label{lem_B}
	Let $$\mathscr{V}=\{v\in C_0^{\infty}(\mathbb{R}^3)~|~\mathrm{div}~v = 0\},~~H=~closure~of~\mathscr{V}~in~L^2(\mathbb{R}^3).$$
	Let $u\in L^1(\mathbb{R}^3)\cap H.$ Then
	$$\int_{\mathbb{R}^3} u dx =0.$$
\end{lemma}
Using the definition of Fourier transform, from Lemma \ref{lem_B} and assumption $u_0\in H^3(\mathbb{R}^3)\cap L^1(\mathbb{R}^3)$ in Theorem 1.2 (ii), one has
$$\hat{u}_0|_{\xi = 0} =\int_{\mathbb{R}^3} u_0 dx = 0,$$
which is contradictory with the assumption $\inf_{0\leq |\xi|\leq R}|\hat{u}_0|\geq c_{0}>0.$ Therefore, the assumption of Theorem 1.2 (ii) in \cite{HWWZ_2022_JDE} was too  restricted.

However, we can fix the above problem  via replacing the assumption ``$(u_0,\tau_0)\in L^1(\mathbb{R}^3)$'' by a weaker one ``$(\hat{u}_0,\hat{\tau}_0)\in L^\infty(\mathbb{R}^3)$''. (For the detailed explanation, see Remark \ref{rem_u_0} right behind Theorem \ref{theorem_decay_m}).  In addition, only minor changes are needed to make the original proof still valid.

\section{The Corrigendum}

To begin with, Theorem 1.2 in \cite{HWWZ_2022_JDE} should be corrected as follows. For the reader’s convenience, the
equation numbers in this section match the corresponding ones in \cite{HWWZ_2022_JDE}.

\begin{theorem}[the corrected Theorem 1.2 of  \cite{HWWZ_2022_JDE}]
	\label{theorem_decay_m}
	For any given $\epsilon$ and $\mu$ satisfying the Case I or Case II, letting $(u^{\epsilon,\mu},\tau^{\epsilon,\mu})$ be the solution as in Theorem 1.1, then the following optimal time-decay estimates hold.

	\medskip

	\noindent {\em (i)} Assume that $ { (\hat{u}_0,\hat{\tau}_0)\in L^\infty(\mathbb{R}^3)}$. Then we have upper time-decay estimates of the solution as below:
	\begin{equation}\label{opti1_m}
	\ \|\nabla^ku^{\epsilon,\mu}(t)\|_{L^2}\leq C_2(1+t)^{-\frac34-\frac{k}{2}},\ k=0,1,2, \tag{1.7}
	\end{equation}
	and
	\begin{equation}\label{opti2_m}
	\ \|\nabla^{k_1}\tau^{\epsilon,\mu}(t)\|_{L^2}\leq C_2(1+t)^{-\frac54-\frac{k_1}{2}},\ k_1=0,1,  \tag{1.8}
	\end{equation}
	for any $t\geq0$ and generic positive constant $C_2$ which depends only on $\|(u_0,\tau_0)\|_{H^3}$, ${ \|(\hat{u}_0,\hat{\tau}_0)\|_{L^\infty_\xi}}$ and $C_1$.
	
\medskip

	\noindent {\em(ii)} Assume that $ { (\hat{u}_0,\hat{\tau}_0)\in L^\infty(\mathbb{R}^3) }$ and in addition that $\inf_{0\leq |\xi|\leq R}|\hat{u}_0|\geq c_{0}>0$, { for some $R>0$}. Then there exist a positive time $t_0=t_0\big(\alpha,\kappa,\beta, {  \|(\hat{u}_0,\hat{\tau}_0)\|_{L^\infty_\xi}}\big)$ and a positive generic constant $c=c(\alpha,\kappa,\beta,c_0,C_2)$, such that
	\begin{equation}\label{opti3_m}
	\|\nabla^ku^{\epsilon,\mu}(t)\|_{L^2}\geq c(1+t)^{-\frac34-\frac{k}{2}},\ k=0,1,2, \tag{1.9}
	\end{equation}
	and
	\begin{equation}\label{opti4_m}
	\|\nabla^{k_1}\tau^{\epsilon,\mu}(t)\|_{L^2}\geq c(1+t)^{-\frac54-\frac{k_1}{2}},\ k_1=0,1, \tag{1.10}
	\end{equation}
	for any $t\geq t_0$.
	
\end{theorem}
\begin{remark}\label{rem_u_0}
	As pointed out by  Schonbek, Schonbek and S\"uli (page 721 in \cite{Schonbek_Suli_1996}), there exists $u_0$ satisfying all the assumptions in Theorem \ref{theorem_decay_m}. Actually, for a given function $g\in C_0^\infty(\mathbb{R}^3)$ with $g(0) =2c_0\neq 0$. Define
	\begin{align}\label{def_u_0}
	\hat{u}_0(\xi)=(\hat{u}_{0,1},\hat{u}_{0,2},\hat{u}_{0,3})^\top:=\left(\frac{\xi_2}{\sqrt{\xi_1^2+\xi_2^2}}g(\xi),-\frac{\xi_1}{\sqrt{\xi_1^2+\xi_2^2}}g(\xi),0\right)^\top.
\end{align}
 Then,
	\begin{align*}
		\|u_0\|_{H^3}^2 =&\, \int_{\mathbb{R}^3}(1+|\xi|^2)^3|\hat{u}_0(\xi)|^2 d\xi = \int_{\mathbb{R}^3}(1+|\xi|^2)^3|g(\xi)|^2 d\xi < \infty, \\
		\widehat{\mathrm{div} u_0} =&\, i\xi_1\hat{u}_{0,1} +  i\xi_2\hat{u}_{0,2} + i\xi_3\hat{u}_{0,3} = i\frac{\xi_1\xi_2}{\sqrt{\xi_1^2+\xi_2^2}}g(\xi)-i\frac{\xi_1\xi_2}{\sqrt{\xi_1^2+\xi_2^2}}g(\xi) = 0,
     \end{align*}
	
	and
	$$\sup_{\xi\in\mathbb{R}^3}|\hat{u}_0| \leq \sup_{\xi\in\mathbb{R}^3}|g(\xi)|<\infty. $$
	Thus, $u_0\in H^3$ with $\sup_{\xi\in\mathbb{R}^3}|\hat{u}_0|<\infty$ and $\mathrm{div}u_0 = 0$.  Moreover, by the continuity of $g$, there exists a constant $R$, such that
	$$
		|g(\xi)| \geq c_0,~~\text{for all }~~0\leq |\xi|\leq R.
	$$
	Hence,
	$$
	\inf_{0\leq |\xi|\leq R}|\hat{u}_0|  = \inf_{0\leq |\xi|\leq R}\left(\frac{\xi_2^2}{\xi_1^2+\xi_2^2}|g(\xi)|^2 + \frac{\xi_1^2}{\xi_1^2+\xi_2^2}|g(\xi)|^2\right)^{\frac{1}{2}} = \inf_{0\leq |\xi|\leq R}|g(\xi)| \geq  c_0 >0.
	$$
	Therefore, $u_0$ defined by \eqref{def_u_0} satisfies all the assumptions  in Theorem \ref{theorem_decay_m}.
\end{remark}
Next, we would like to remark that the proof of Theorem 1.2 in \cite{HWWZ_2022_JDE} can be corrected via replacing the quantities quoted in the second column by the corresponding ones in the third column in the following Table.

\begin{table}[h!]
	\begin{center}
		\caption{Corrigendum}
	\begin{tabular}{|c|c|c|}
		\hline
		\textbf{Page}& \textbf{Original form }&  \textbf{Corrected form} \\
		\hline
		\multirow{2}{*}{475} & ``${  \| {U}_0\|_{L^p}}$'' in (4.5), (4.7) and the estimate right before&  \multirow{2}{*}{${  \|\hat{U}_0\|_{L^q_\xi}}$}\\
		&     (4.7).& \\
		\hline
		\multirow{2}{*}{477--478} & \multirow{2}{*}{``${ \| {U}_0\|_{L^1}}$'' in (4.11) and (4.13).}  & \multirow{2}{*}{${  \|\hat{U}_0\|_{L^\infty_\xi}}$}\\
		&   &\\
		\hline
		 {478--479,} &  ``${  \| {U}_0\|_{L^1}^2}$'' in (4.12), (4,19), (5.4), the estimates right before     & \multirow{2}{*}{${  \|\hat{U}_0\|_{L^\infty_\xi}^2}$}\\
		487 & (4.13) and (4.17). &\\
		\hline
		\multirow{3}{*}{482, 484} &the line under (4.31) and the first line in page 484, ``${ \|(u_0,\tau_0)\|_{L^1}}$'' in      & \multirow{3}{*}{${ \|(\hat{u}_0,\hat{\tau}_0)\|_{L^\infty_\xi}}$}\\
		&$t_2 = t_2(\alpha,\kappa,\beta,{ \|(u_0,\tau_0)\|_{L^1}})$, $c_2 = c_2(\alpha,\kappa,\beta,c_0,{ \|(u_0,\tau_0)\|_{L^1}})$ &\\
		&and $t_0 = t_0(\alpha,\kappa,\beta,{ \|(u_0,\tau_0)\|_{L^1}}).$ &\\
		\hline
		 {482--483} &  ``${ \|(u_0,\sigma_0)\|_{L^1}}$'' in (4.34) and  the estimate right behind (4.33).&  {${ \|(\hat{u}_0,\hat{\tau}_0)\|_{L^\infty_\xi}}$}\\
		\hline
	\end{tabular}
\end{center}
\end{table}


\end{document}